\documentclass[12pt]{article}
\usepackage{amsfonts}
\usepackage{amsmath}
\usepackage{amssymb}
\usepackage{amsthm}

\def\Z{\Bbb Z}

\def\C{\Bbb C}

\def\l{\left}
\def\r{\right}
\def\bg{\bigg}
\def\({\bg(}
\def\){\bg)}
\def\]{]\!]}
\def\[{[\![}
\def\t{\text}

\def\em{\emptyset}

\def\ov{\overline}
\def\ls{\leqslant}
\def\gs{\geqslant}
\def\se {\subseteq}

\def\ra{\rightarrow}

\def\supp {{\rm supp}}

\theoremstyle{plain}
\newtheorem{theorem}{Theorem}

\theoremstyle{definition}

\def\Ack{\medskip\noindent {\bf Acknowledgment.}}
\theoremstyle{remark} 

\begin{document}
\hbox{J. Number Theory 129(2009), no.\,2, 434--438}
\bigskip
\centerline {\bf A Variant of Tao's Method with Application to Restricted Sumsets}
\bigskip
\centerline{Song Guo$^1$ and Zhi-Wei Sun$^2$}
\medskip
\centerline{$^1$Department of Mathematics, Huaiyin Teachers College}
\centerline{Huaian 223300, People's Republic of China}
\centerline{guosong77@hytc.edu.cn}
\medskip
\centerline{$^2$Department of Mathematics, Nanjing University}
\centerline{Nanjing 210093, People's Republic of China}
\centerline{zwsun@nju.edu.cn}
\centerline{\tt http://math.nju.edu.cn/$\sim$zwsun}
\medskip
\abstract {In this paper, we develop Terence Tao's harmonic analysis method
and apply it to restricted sumsets.  The well known Cauchy-Davenport
theorem asserts that if $\em\neq A,B\se \Z/p\Z$ with $p$ a prime,
then $|A+B|\gs \min\{p,\,|A|+|B|-1\}$, where $A+B=\{a+b:\,a\in A,\,b\in
B\}$. In 2005, Terence Tao gave a harmonic analysis proof of the
Cauchy-Davenport theorem, by applying a new form of the uncertainty
principle on Fourier transform. We modify Tao's
method so that it can be used to prove the following extension of the
Erd\H{o}s-Heilbronn conjecture: If $A,B,S$
are nonempty subsets of $\Z/p\Z$ with $p$ a prime, then $\big|\{a+b:\,a\in
A,\,\,b\in B,\,a-b\not\in S\}\big|\gs \min \{p,\,|A|+|B|-2|S|-1\}$.}

\footnote {
\newline\indent\ \ The second author is supported by the National Natural Science
Foundation of People's Republic of China.}

 \vspace{4mm}

\section{Introduction}

 Let $p$ be a prime, and let $A$ and $B$ be two subsets of the finite field
 $$\Z_p=\Z/p\Z=\{\bar r=r+p\Z:\ r\in\Z\}.$$
 Set
\begin{equation}A+B=\{a+b:\ a\in A,\, b\in B\}
\end{equation}and
\begin{equation}A\dot{+}B=\{a+b:\ a\in A,\, b\in B,\, a\neq b\}.
\end{equation}
The well-known Cauchy-Davenport theorem asserts that
\begin{equation} |A+B|\gs \min\{p,\,|A|+|B|-1\}.\end{equation}
In 1964 P. Erd\H{o}s and H. Heilbronn \cite{EH} conjectured that
\begin{equation}|A\dot{+}A|\gs \min\{p,\,2|A|-3\};\end{equation}
this was confirmed by J. A. Dias da Silva and Y. O. Hamidoune
\cite{SH} in 1994. In 1995-1996 N. Alon, M. B. Nathanson and
I. Z. Ruzsa \cite{ANR} proposed the so-called polynomial method to handle similar
problems. By the powerful polynomial
method (cf. \cite{A1} and \cite{ANR}),
many interesting results on restricted sumsets have been obtained (see, e.g., \cite{HSu},
\cite{LS}, \cite{PSu}, \cite{PS}, \cite{Sun}).

In 2005, Terence Tao \cite{T} developed a harmonic analysis method in this
area, applying a new form of the uncertainty principle on Fourier
transform. Let $p$ be a prime. For a complex-valued
function $f: \Z_p\rightarrow \C$, we define its support $\supp (f)$ and its Fourier
transform $\hat{f}:\Z_p\ra \C$ as follows:
\begin{equation}
\supp (f)=\{x\in \Z_p: \ f(x)\neq 0\}\end{equation} and
\begin{equation}
\hat{f}(x)=\sum_{a\in \Z_p}f(a)e_p(ax)\quad\t{for all}\ x\in \Z_p,\end{equation}
where $e_p(\bar r)=e^{-2\pi ir/p}$ for any $r\in\Z$.

Here is the main result of the paper \cite{T}.

\begin{theorem} {\rm (T. Tao \cite{T})} Let $p$ be an odd prime. If $f : \Z_p \ra
\C$ is not identically zero, then
\begin{equation}|\supp (f)|+|\supp (\hat{f})|\gs p+1.
\end{equation}
Moreover, given two non-empty subsets A and B of $\Z_p$ with $|A|+|B|\gs p+1$,
we can find a function $f : \Z_p \ra \C$ such that $\supp (f)=A$ and
$\supp (\hat{f})=B$.
\end{theorem}

Using this theorem Tao \cite{T} gave a new proof of Cauchy-Davenport theorem.
Note that the inequality (7) was also discovered independently by
Andr\'{a}as Bir\'{o} (cf. \cite{Fr} and \cite{T}). In
this article we adapt the method further and use the refined method to deduce
the following result.

\begin{theorem} Let $A$ and $B$ be non-empty subsets of $\Z_p$ with $p$ a prime,
and let
\begin{equation}C=\{a+b:\ a\in A,\, b\in B,\, a-b\not\in S\}
\end{equation} with $S\se \Z_p$.  Then we have
\begin{equation}|C|\gs \min \{p,\,|A|+|B|-2|S|-1\}.\end{equation}
\end{theorem}

Theorem 2 in the case $S=\em$ reduces to the Cauchy-Davenport theorem.
When $A=B$ and $S=\{0\}$, Theorem 2 yields the Erd\H{o}s-Heilbronn conjecture.
In the case $p\not=2$ and $\em\not=S\subset\Z_p$, Pan and Sun \cite[Corollary 2]{PSu}
obtained the stronger inequality
\begin{equation}|C|\gs \min \{p,\,|A|+|B|-|S|-2\}\end{equation}
via the polynomial method. The second author (cf. \cite{HSu}) ever conjectured that 2 in (10)
can be replaced by 1 if $|S|$ is even.
We conjecture that when $A\not=B$ we can also substitute 1 for 2 in (10).

\medskip
\section{Proof of Theorem 2}
Without loss of generality, we let  $|A|\ls |B|$. When $|A|+|B|\ls
2|S|+1$ or $|A|=1$, (9) holds trivially. Below we suppose that
$|A|+|B|> 2|S|+1$ and $|A|\gs 2$.

 In the case $p=2$, we have $A=B=\Z_2$ and $C=(A+B)\setminus S=\Z_2\setminus S$, thus
 $$|C|=2-|S|\gs\min\{2,\,|A|+|B|-2|S|-1\}=\min\{2,\, 3-2|S|\}.$$

Below we assume that $p$ is an odd prime.
Set $k=p-|A|+1\in [1,p-1]$ and
$l=p-|B|+1\in [1,p-1]$. Then $k+l\ls 2p-2|S|$ and $l\ls p-|S|$ since
$2|B|\gs |A|+|B|\gs 2|S|+2$.
Define
\begin{equation}\hat{A}=\{\bar{0},\ldots,\ov{k-1}\}=\{\bar0,\ldots,\ov{p-|A|}\}
\end{equation}
and
\begin{equation}\hat{B}=\{\ov{p-|S|-l+1},\ldots,\ov{p-|S|}\}
=\{\ov{|B|-|S|},\ldots,\ov{p-|S|}\}.
\end{equation}
Clearly, $|\hat{A}|=p+1-|A|$ and $|\hat{B}|=p+1-|B|$. By Theorem 1
there are functions $f,g:\Z_p\ra \C$ such that
\begin{equation}\supp(f)=A,\ \supp (\hat{f})=\hat{A},\ \supp (g)=B,\ \supp(\hat{g})=\hat{B}.
\end{equation}

 Now we define a function $F:\Z_p\ra \C$ by
\begin{equation}F(x)=\sum_{a\in \Z_p}f(a)g(x-a)
\prod_{d\in S}(e_p(x-a)-e_p(a-d)).
\end{equation}
For each $x\in \supp (F)$, there exists $a\in\supp(f)$ with $x-a\in \supp (g)$
and $d:=a-(x-a)\not\in S$, hence $x=a+(x-a)\in C$. Therefore
\begin{equation}\supp(F)\se C.\end{equation}

For any $x\in\Z$ we have
$$\hat{F}(x)=\sum_{b\in\Z_p}F(b)e_p(bx)
 =\sum_{a\in \Z_p}\sum_{b\in
 \Z_p}f(a)g(b-a)e_p(bx)P(a,b),$$
where
\begin{align*}P(a,b)=&\prod_{d\in S}(e_p(b-a)-e_p(a-d))
\\=&\sum_{T\se S}(-1)^{|T|}
 e_p\l((|S|-|T|)(b-a)\r)e_p\bigg(|T|a-\sum_{d\in T}d\bigg).
 \end{align*}
Therefore
\begin{align*} \hat{F}(x)
 =&\sum_{T\se S}(-1)^{|T|}e_p\bigg(-\sum_{d\in T}d\bigg)
 \sum_{a\in \Z_p}f(a)e_p(ax+|T|a) \\
 &\quad\times\sum_{b\in\Z_p}g(b-a)e_p\l((b-a)x+(|S|-|T|)(b-a)\r) \\
=&\sum_{T\se S}(-1)^{|T|}e_p\bigg(-\sum_{d\in
T}d\bigg)\hat{f}\l(x+\ov{|T|}\r)\hat{g}\l(x+\ov{|S|-|T|}\r) .
\end{align*}
For $T\se S$, if $\ov{p-|S|}+\ov{|S|-|T|}\in \supp(\hat{g})=\hat{B}$,
then we must have $|T|=|S|$ (i.e., $T=S$) by the definition of $\hat B$.
It follows that
 $$\hat{F}\l(\ov{p-|S|}\r)=(-1)^{|S|}e_p\bigg(-\sum_{d\in S}d\bigg)
 \hat{f}\l(\ov{0}\r)\hat{g}\l(\ov{p-|S|}\r)\neq 0$$
since $\ov{0}\in\hat{A}=\supp (\hat{f})$
and $\ov{p-|S|}\in\hat{B}=\supp (\hat{g})$. With the helps of
(15) and Theorem 1, we get
$$|C|\gs |\supp (F)|\gs  p+1-|\supp (\hat{F})|.$$

Suppose that $x\in \supp (\hat{F})$. By the above, there is a subset $T$
of $S$ such that $x+\ov{|T|}\in \supp (\hat{f})=\hat A$ and
$x+\ov{|S|-|T|}\in \supp(\hat{g})=\hat B$. As $0\ls
|T|\ls|S|$,
$$x+\ov{|T|}\in \hat{A} \Longrightarrow x\in
\{\ov{p-|S|},\ldots,\ov{p-1},\ov0,\ldots,\ov{k-1}\}$$
and
$$x+\ov{|S|-|T|}\in \hat{B}
\Longrightarrow x\in \{\ov{|B|-2|S|},\ldots,\ov{p-|S|}\}.$$
Therefore $x=\ov{p-|S|}$, or $x=\ov{r}$ for some $r\in [|B|-2|S|,\,k-1]$.

 If $|A|+|B|\gs p+2|S|+1$, then $k-1=p-|A|<|B|-2|S|$,
hence $\supp (\hat{F})=\{\ov{p-|S|}\}$ and thus $|C|\gs p$. If
$|A|+|B|< p+2|S|+1$, then
$$|\supp (\hat{F})|\ls 1+k-(|B|-2|S|)= k+l-p+2|S|$$
and hence
$$|C|\gs p+1-k-l+p-2|S|=|A|+|B|-2|S|-1.$$
So (9) always holds. We are done.

\smallskip

\Ack\  The authors are grateful to the referee for his/her helpful
comments.

\end{document}